\documentclass[12pt]{article}
\usepackage{latexsym,amsfonts}

\newtheorem{thm}{Theorem}[section]
\newtheorem{prop}[thm]{Proposition}
\newtheorem{lem}[thm]{Lemma}

\catcode`\@=11
\@addtoreset{equation}{section}
\catcode`\@=12

\newcommand\goa{{\mathfrak a}}
\newcommand\gog{{\mathfrak g}}
\newcommand\goh{{\mathfrak h}}
\newcommand\gok{{\mathfrak k}}
\newcommand\gol{{\mathfrak l}}
\newcommand\gom{{\mathfrak m}}
\newcommand\gon{{\mathfrak n}}
\newcommand\gop{{\mathfrak p}}
\newcommand\gor{{\mathfrak r}}
\newcommand\gou{{\mathfrak u}}

\newcommand\DD{{\mathbb D}}
\newcommand\RR{{\mathbb R}}
\newcommand\CC{{\mathbb C}}

\newcommand\Proof{{\bf Proof}\quad}
\newcommand\bLP{\bigbreak\noindent}

\newcommand\ad{{\rm ad}}
\newcommand\Ad{{\rm Ad}}
\newcommand\supp{{\rm supp}}

\font\addressfont=cmr10 at 10truept
\font\sladdressfont=cmsl10 at 10truept
\font\ttaddressfont=cmtt10 at 10truept

\begin{document}

\title {Invariant Differential Operators on Nonreductive Homogeneous
Spaces}
\author{Tom H. Koornwinder}
\date{}

\maketitle

\begin{abstract}\noindent
A systematic exposition is given of the theory of invariant differential
operators on a not necessarily reductive homogeneous space.
This exposition is modelled on Helgason's treatment of the general
reductive case and the special nonreductive case of the space of
horocycles.
As a final application the differential operators on (not a priori
reductive) isotropic pseudo-Riemannian spaces are characterized.
\end{abstract}

\bLP

\bLP

\bLP
{\sl Note\/}:
This is an essentially unchanged electronic version of
Report ZW 153/81, Mathematisch Centrum, Amsterdam, 1981;
MR 82g:43011.

\section{Introduction}
Let $G$ be a Lie group and $H$ a closed subgroup.
Let $\gog$ and $\goh$ denote the corresponding Lie algebras.
Suppose that the coset space $G/H$ is {\sl reductive}, i.e., there is a
complementary subspace $\gom$ to $\goh$ in $\gog$ such that
$\Ad_G(H)\gom\subset \gom$.
Let $\DD(G/H)$ denote the algebra of $G$-invariant differential
operators on $G/H$.
The main facts about $\DD(G/H)$ are summarized below (cf.\ HELGASON
[3, Ch.III], [4, Cor. X.2.6, Theor. X.2.7], [6, \S 2]).

Let $\DD(G)$ be the algebra of left invariant differential operators
on $G$, $\DD_H(G)$ the subalgebra of operators which are right
invariant under $H$ and $S(\gog)$ the complexified symmetric algebra
over $\gog$.
Let $\lambda\colon S(\gog)\to\DD(G)$ denote the symmetrization
mapping.
$I(\gom)$ denotes the set of $\Ad_G(H)$-invariants in $S(\gom)$.
Then
\begin{equation}
\DD_H(G)=\DD(G)\goh\cap \DD_H(G)\oplus\lambda(I(\gom)).
\end{equation}
Let $\pi\colon G\to G/H$ be the natural mapping.
Let $C^\infty_H(G)$ consist of the $C^\infty$-functions on $G$ which are
right invariant under $H$.
Write ${\tilde f}:=f\circ\pi$ ($f\in C^\infty(G/H)$) and
$(D_uf)^{\sim} :=u{\tilde f}$ ($f\in C^\infty(G/H)$, $u\in\DD_H(G)$).
Then $D_u\in\DD(G/H)$.
\begin{thm}\quad
The mapping $u\mapsto D_u$ is an algebra homomorphism from $\DD_H(G)$
onto $\DD(G/H)$ with kernel $\DD(G)\goh\cap\DD_H(G)$.
The mapping $P\mapsto D_{\lambda(P)}\colon I(\gom)\to\DD(G/H)$ is a linear
bijection.
\end{thm}

\medskip
Theorem 1.1. is of basic importance for the analysis on symmetric spaces.
In particular, it can be shown that $\DD(G/H)$ is commutative if
$G/H$ is a pseudo-Riemannian symmetric space which admits a relatively
invariant measure.
In its most general form this result was proved by DUFLO [1] in an
algebraic way.
G. van Dijk kindly communicated a short analytic proof of Duflo's result
to me (unpublished).
In [1] DUFLO used generalizations of (1.1) and Theorem 1.1 to the case of
homogeneous line bundles over $G/H$.
These can be proved by only minor changes of Helgason's original proofs.

There exist nonreductive coset spaces $G/H$ for which $\DD(G/H)$ is
still commutative.
For instance, let $G$ be a connected real semisimple Lie group and let
$M$ and $N$ be the usual subgroups of $G$.
Then $G/MN$ is the space of horocycles and $\DD(G/MN)$ is commutative.
In order to prove this, formula (1.1) and Theorem 1.1 have to be adapted
to the nonreductive case.
While HELGASON [5, \S 4], [6, \S 3] has done this in an ad hoc way for the
special coset spaces under consideration, it is the purpose of the
present note to give a more systematic exposition of the theory of
$\DD(G/H)$ for a not necessarily reductive coset space.

Furthermore, following Duflo, the theory will be developed for invariant
differential operators on homogeneous line bundles over $G/H$.
As a final application we will characterize $\DD(G/H)$ for isotropic
pseudo-Riemannian symmetric spaces $G/H$ without a priori knowledge that
$G/H$ is reductive.
Throughout HELGASON [4] will be our standard reference.

\section{Development of the general theory}
Let $G$ be a Lie group with Lie algebra $\gog$.
For $X\in \gog$ define the vector field ${\tilde X}$ on $G$ by
\begin{equation}
({\tilde X}f)(g):=\frac{d}{dt}f(g\exp tX)\big|_{t=0}\,, \quad f\in
C^\infty(G), \ g\in G.
\end{equation}
Then the mapping $X\mapsto {\tilde X}$ is an isomorphism from $\gog$ onto the
Lie algebra of left invariant vector fields on $G$.
Throughout this section let $X_1,\dots,X_n$ be a fixed basis of $\gog$.

For a finite-dimensional real vector space $V$ the symmetric algebra
$S(V)$ is defined as the algebra of all polynomials with complex
coefficients on $V^*$, the dual of $V$.
Let $S^m(V)$ respectively $S_m(V)$ ($m=0,1,2,\dots\;$) denote the space of
homogeneous polynomials of degree $m$ on $V^*$, respectively of
polynomials of degree $\le m$ on $V^*$.
Thus $S^m(G)$ is spanned by the monomials $X_{i_1}X_{i_2}\dots X_{i_m}$
($i_1,\dots,i_m\in \{1,\dots, n\}$).

Let $\DD(G)$  be the algebra of left invariant differential operators
on $G$ with complex coefficients.
For $P\in S(\gog)$ define an operator $\lambda(P)$ on $C^\infty(G)$ by
\begin{equation}
(\lambda(P)f)(g):=P\left(\frac{\partial}{\partial
t_1},\cdots,\frac{\partial}{\partial
t_n}\right)f(g\exp(t_1X_1+\cdots+t_nX_n))\Big|_{t_1=\cdots=t_n=0},
\end{equation}
where
\[
P\left(\frac{\partial}{\partial t_1},\cdots,\frac{\partial}{\partial
t_n}\right):=\frac{\partial^m}{\partial t_{i_1}\cdots \partial t_{i_m}}\quad
\mbox{for $P=X_{i_1}\cdots X_{i_m}$.}
\]
It is proved in [4, Prop. II.1.9 and p. 392] that:
\begin{prop}
The mapping $P\mapsto\lambda(P)$ is a linear bijection from $S(\gog)$ onto
$\DD(G)$. It satisfies
\begin{equation}
\lambda(Y^m)={\tilde Y}^m, \quad Y\in \gog;
\end{equation}
\begin{equation}
\lambda(Y_1\dots Y_m)=\frac{1}{m!}\sum_{\sigma\in S_m} {\tilde
Y}_{\sigma(1)}\dots{\tilde Y}_{\sigma(m)}, \quad Y_1,\dots, Y_m\in \gog.
\end{equation}
The definition of $\lambda$ is independent of the choice of the basis of
$\gog$.
\end{prop}

\medskip
The mapping $\lambda$  is called {\sl symmetrization}.
The Lie algebra $\gog$ is embedded  as a subspace of $\DD(G)$ under the
mapping $X\to{\tilde X}$.
Any homomorphism from $\gog$ to $\gog$ uniquely extends to a homomorphism from
$\DD(G)$ to $\DD(G)$ and any linear mapping from $\gog$ to $\gog$
uniquely extends to a homorphism from $S(\gog)$ to $S(\gog)$.
In particular, for $g\in G$, the automorphism $\Ad(g)$  of $\gog$ uniquely
extends to automorphisms of both $S(\gog)$ and $\DD(G)$ and
\begin{equation}
\lambda(\Ad(g)P)=\Ad(g)\lambda(P), \quad P\in S(\gog), \ g\in G.
\end{equation}
For $g,g_1\in G, \ f\in C^\infty(G), \ D\in\DD(G)$ write
\[f^{R(g)}(g_1):=f(g_1g); \quad D^{R(g)}f:=(Df^{R(g^{-1})})^{R(g)}.\]
Then
\begin{equation}
\Ad(g)D=D^{R(g^{-1})}, \quad D\in\DD(G), \ g\in G.
\end{equation}

Let $H$ be a closed subgroup of $G$ and let $\goh$ be the corresponding
subalgebra.
Let $\gom$ be a subspace of $\gog$ complementary to $\goh$.
Let $X_1,\dots, X_r$ be a basis of $\gom$ and $X_{r+1},\dots, X_n$ a basis
of $\goh$.
Let $\chi$ be a character of $H$, i.e.\ a continuous homomorphism from
$H$ to the multiplicative group $\CC{\backslash\{0\}}$.
Throughout this section, $H$, $\gom$, the basis and $\chi$ will be assumed
fixed.

Let $\pi\colon G\to G/H$ be the canonical mapping.
Write $0:=\pi(e)$.
Let
\begin{equation}
C^\infty_{H,\chi}(G):=\{f\in C^\infty(G)\mid f(gh)=f(g)\chi(h^{-1}), \
g\in G, \ h\in H\}.
\end{equation}

Sometimes we will assume that $\chi$ has an extension to a character on
$G$.
This assumption clearly holds if $\chi\equiv 1$ on $H$, but it does not
hold for general $H$ and $\chi$.
For instance, if $G=SU(2)$ or $SL(2,\RR)$ and $H=SO(2)$ then
nontrivial characters on $H$ do not extend to characters on $G$.

If $\chi$ extends to a character on $G$ then we define a linear
bijection $f\mapsto{\tilde f}\colon C^\infty(G/H)\to C^\infty_{H,\chi}(G)$ by
\begin{equation}
{\tilde f}(g):=f(\pi(g))\chi(g^{-1}), \quad g\in G.
\end{equation}
\begin{lem}
Let $P\in S(m)$.
If $\lambda(P)f=0$ for all $f\in C^\infty_{H,\chi}(G)$ then $P=0$.
\end{lem}
\medskip
\noindent
\Proof
For each $f\in C^\infty(G/H)$ we can find $F\in
C^\infty_{H,\chi}(G)$ such that
\[F(\exp(t_1X_1+\cdots+t_rX_r))=f(\exp(t_1X_1+\cdots+t_rX_r)\cdot 0)\]
for $(t_1,\dots,t_r)$ in some neighbourhood of $(0,\dots,0)$.
Hence
\[
0=(\lambda(P)F)(e)
=P\left(\frac{\partial}{\partial t_1},\cdots,\frac{\partial}{\partial
t_r}\right)f(\exp(t_1X_1+\cdots+t_rX_r)\cdot 0\Big|_{t_1=\cdots=t_r=0}
\]
for all $f\in C^\infty(G/H)$, so $P=0$.\hfill $\Box$

\bigskip
Let the differential of $\chi$ also be denoted by $\chi$.
Let $\goh^\CC$ be the complexification of $\goh$.
Let
\begin{equation}
\goh^\chi:=\{X+\chi(X)\mid X\in \goh^\CC\}\subset\DD(G).
\end{equation}
Clearly, $Df=0$ if $f\in C^\infty_{H,\chi}(G)$ and $D\in \goh^\chi$.
Let $\DD(G)\goh^\chi$ be the linear span of all $vw$ with $v\in\DD(G),
\ w\in \goh^\chi$.
Observe that, by Proposition 2.1, ${\tilde Y}_1\dots{\tilde
Y}_m\in\lambda(S_m(\gog))$ for $Y_1,\dots, Y_m\in \gog$.
The following proposition was proved in [4, Lemma X.2.5] for
$\chi\equiv 1$.
\begin{prop}
There are the direct sum decompositions
\begin{equation}
\lambda(S_m(\gog))=\lambda(S_{m-1}(\gog))\goh^\chi\oplus\lambda(S_m(\gom))
\end{equation}
and
\begin{equation}
\DD(G)=\DD(G)\goh^\chi\oplus\lambda(S(\gom)).
\end{equation}
\end{prop}

\medskip
\noindent
\Proof 
First we prove by complete induction with respect to $\gom$ that
\[\lambda(S_m(\gog))\subset\lambda(S_{m-1}(\gog))\goh^\chi+
\lambda(S_m(\gom)).\]
This clearly holds for $m=0$.
Suppose it is true for $m<d$.
Let
\[P=X_1^{d_1}\cdots X_n^{d_n}, \quad d_1+\cdots+d_n=d.\]
If $d_{r+1}\cdots+d_n=0$, then $P\in S_d(\gom)$, so
$\lambda(P)\in\lambda(S_d(\gom))$.
If $d_{r+1}+\cdots+d_n>0$ then, by (2.4), $\lambda(P)$ is a linear
combination of certain elements ${\tilde Y}_1\cdots{\tilde Y}_d$ with
$Y_i\in h$ for at least one $i$, so
\[\lambda(P)\in\lambda(S_{d-1}(\gog))\goh^\CC+
\lambda(S_{d-1}(\gog))\subset\lambda(S_{d-1}(\gog))\goh^\chi+
\lambda(S_{d-1}(\gog)).\]
Now apply the induction hypothesis.
This yields (2.10) and (2.11) (use Proposition 2.1) except for the
directness.

To prove the directness of the sum (2.11), suppose that $P\in S(\gom)$ and
$\lambda(P)\in \DD(G)\goh^\chi$.
Then $\lambda(P)f=0$ for all $f\in C^\infty_{H,\chi}(G)$, so $P=0$ by
Lemma 2.2. \hfill $\Box$
\begin{lem}
Let $D\in \DD(G)$.
Then $Df=0$ for all $f\in C^\infty_{H,\chi}(G)$ if and only if
$D\in\DD(G)\goh^\chi$.
\end{lem}

\medskip\noindent
\Proof 
Apply Proposition 2.3 and Lemma 2.2. \hfill $\Box$

\bigskip
Let us define
\begin{equation}
\DD_{H,\chi,{\rm mod}}(G):=\{D\in\DD(G)\mid \Ad(h)D-D\in\DD(G)\goh^\chi \ \mbox{for all} \ h\in H\}.
\end{equation}
This definition is motivated by the following lemma.
\begin{lem}
Let $D\in\DD(G)$.
Then the following two statements are equivalent.
\begin{itemize}
\item[{\rm(i)}] $D\in\DD_{H,\chi,{\rm mod}}(G)$.
\item[{\rm(ii)}] $f\in C^\infty_{H,\chi}(G)\Rightarrow Df\in
C^\infty_{H,\chi}(G)$.
\end{itemize}
\end{lem}
\medskip\noindent
\Proof 
Let $D\in \DD(G)$.
If  $f\in C^\infty_{H,\chi}(G), \ h\in H$ then
\begin{itemize}
\item [($\star$)]\qquad\qquad
$(Df)^{R(h)}=D^{R(h)}f^{R(h)}=\chi(h^{-1})D^{R(h)}f$.
\end{itemize}

\medskip\noindent
First assume (i).
If $f\in C^\infty_{H,\chi}(G), \ h\in H$, then
$(D^{R(h)}-D)f=(\Ad(h)D-D)f=0$, so combination with ($\star$) yields
$(Df)^{R(h)}=\chi(h^{-1})Df$, i.e., $Df\in C^\infty_{H,\chi}(G)$.
Conversely, assume (ii).
If $f\in C^\infty_{H,\chi}(G), \ h\in H$, then
$(Df)^{R(h)}=\chi(h^{-1})Df$, so combination with ($\star$) yields
$(D^{R(h)}-D)f=0$.
Hence $\Ad(h)D-D=D^{R(h)}-D\in \DD(G)\goh^\chi$ by Lemma 2.4.\hfill $\Box$

\bigskip
{}From the preceding results the following theorem is now obvious.
\begin{thm}
\quad\par
\begin{itemize}
\item [\rm(a)]  $\DD_{H,\chi,{\rm mod}}(G)$ is a subalgebra of $\DD(G)$.
\item [\rm(b)]  $\DD(G)\goh^\chi$ is a two-sided ideal in
$\DD_{H,\chi,{\rm mod}}(G)$.
\item [\rm(c)]  There is the direct sum decomposition.
\begin{equation}
\DD_{H,\chi,{\rm mod}}(G)=
\DD(G)\goh^\chi\oplus\lambda(S(\gom))\cap \DD_{H,\chi,{\rm
mod}}(G).
\end{equation}
\item [\rm(d)]  Define the mappings $A$ and $B$ by
\begin{eqnarray*}
u\stackrel{A}{\longmapsto} u({\rm mod}\;
\DD(G)\goh^\chi)&\stackrel{B}{\longmapsto}&
u\big|_{C^\infty_{H,\chi}(G)}\colon\\
\lambda(S(\gom))\cap\DD_{H,\chi,{\rm mod}}(G)&\stackrel{A}{\longrightarrow}&
\DD_{H,\chi,{\rm mod}}(G)/\DD(G)\goh^\chi
\stackrel{B}{\longrightarrow}\DD_{H,\chi,{\rm mod}}
\Big|_{C^\infty_{H,\chi}(G)}.
\end{eqnarray*}
Then $A$ is a linear bijection and $B$ is an algebra isomorphism onto.
\end{itemize}
\end{thm}

Define the mapping $\sigma\colon \gog\to \gom$ by
\begin{equation}
\sigma(X+Y):=X, \quad X\in \gom, \ Y\in \goh.
\end{equation}
Consider $S(\gom)$ as a subalgebra of $S(\gog)$.
Thus, if $P\in S(\gom)$ and $h\in H$, then $\Ad(h)P\in S(\gog)$ and
$\sigma\circ
\Ad(h)P\in S(\gom)$ are well-defined.
By an application of (2.4) we see that, if $Q\in S_m(\gog)$, then
\begin{equation}
\lambda(\sigma Q-Q)\in\lambda(S_{m-1}(g))+\DD(G)\goh^\chi.
\end{equation}
Define the algebra
\begin{equation}
I_{\rm mod}(\gom):=
\{P\in S(\gom)\mid\sigma\circ \Ad(h)P=P \ \mbox{for all} \ h\in H\}.
\end{equation}
\begin{lem}
Let $P\in S(\gom)$ such that $\lambda(P)\in\DD_{H,\chi,{\rm mod}}(G)$.
Write $P=P^m+P_{m-1}$, where $P^m\in S^m(\gom), \ P_{m-1}\in S_{m-1}(\gom)$.
Then $P^m\in I_{{\rm mod}}(\gom)$.
\end{lem}

\medskip\noindent
\Proof 
$\lambda(\Ad(h)P-P)\in\DD(G)\goh^\chi$ by (2.12).
Hence
\[\lambda(\Ad(h)P^m-P^m)\in\lambda(S_{m-1}(\gog))+\DD(G)\goh^\chi.\]
So
\[\lambda(\sigma\circ \Ad(h)P^m-P^m)\in\lambda(S_{m-1}(\gog))+
\DD(G)\goh^\chi\subset\lambda(S_{m-1}(\gom))+\DD(G)\goh^\chi,\]
where we used (2.16) and (2.10).
By directness of the decomposition (2.10):
\[\sigma\circ \Ad(h)P^m-P^m\in S_{m-1}(\gom).\]
Hence $\sigma\circ \Ad(h)P^m-P^m$, being homogeneous of degree $m$, is the zero
polynomial. \hfill $\Box$
\begin{prop}
If $\lambda(I_{\rm mod}(\gom))\subset\DD_{H,\chi,{\rm mod}}(G)$
then
\[\lambda(I_{\rm mod}(\gom))=\lambda(S(\gom))\cap\DD_{H,\chi,{\rm mod}}(G)\]
and the mapping
\[D\mapsto D\Big|_{C^\infty_{H,\chi}(G)}\colon\lambda(I_{{\rm mod}}(\gom))\to
\DD_{H,\chi,{\rm mod}}(G)\Big|_{C^\infty_{H,\chi}(G)}\]
is a linear bijection.
\end{prop}
\medskip\noindent
\Proof 
Use complete induction with respect to the degree of $P\in S(\gom)$ in order to
prove that $P\in I_{\rm mod}(\gom)$ if $\lambda(P)\in\DD_{H,\chi,{\rm mod}}(G)$
(apply Lemma 2.7).
The second implication in the proposition follows from Theorem 2.6(d).\hfill
$\Box$

\bigskip
Suppose for the moment that $\chi$ extends to a character on $\gog$ and
remember the
mapping $f\to{\tilde f}$ defined by (2.8).
For $u\in\DD_{H,\chi,{\rm mod}}(G)$ define an operator  $D_u$ acting on
$C^\infty(G/H)$ by
\begin{equation}
(D_uf)^\sim:=u{\tilde f}, \quad f\in C^\infty(G/H).
\end{equation}
Then $\supp(D_uf)\subset\supp(f)$, hence, by Peetre's theorem (cf.\ for instance
NARASIMHAN [7, \S 3.3]), $D_u$ is a differential operator on  $G/H$.
One easily shows that $D_u\in\DD(G/H)$,  the space of $G$-invariant
differential operators on $G/H$.
\begin{thm}
Suppose that $\chi$ extends to a character on $G$.
Then the mapping
\[u\Big|_{C^\infty_{H,\chi}(G)}\stackrel{C}{\longmapsto}D_u\colon\quad
\DD_{H,\chi,{\rm
mod}}(G)\Big|_{C^\infty_{H,\chi}(G)}\stackrel{C}{\longrightarrow}\DD(G/H)\]
is an algebraic isomorphism onto.
\end{thm}
\medskip\noindent
\Proof 
Clearly, $C$ is an isomorphism into.
In order to prove the surjectivity let $D\in\DD(G/H)$.
Then there is a polynomial $P\in S(\gom)$ such that
\[(Df)(g\cdot 0)=
P\left(\frac{\partial}{\partial t_1},\cdots,\frac{\partial}{\partial
t_r}\right)f(g\exp(t_1X_1+\cdots+t_rX_r)\cdot 0\Big|_{t_1=\cdots=t_r=0}\]
for all $f\in C^\infty(G/H)$ and for $g=e$.
By the $G$-invariance of $D$ this formula holds for all $g\in G$.
By (2.8) and (2.2) this becomes
\[\chi(Df)^\sim=\lambda(P)(\chi{\tilde f}), \quad\mbox{i.e.,}\quad
(Df)^\sim=(\chi^{-1}\lambda(P)\circ\chi)({\tilde f}).\]
Clearly, $\chi^{-1}\lambda(P)\circ\chi\in\DD(G)$ and, by Lemma 2.5, we have
$\chi^{-1}\lambda(P)\circ\chi\in\DD_{H,\chi{\rm mod}}(G)$.
Thus, by (2.17),
$D=D_{\chi^{-1}\lambda(P)\circ\chi}$.\hfill $\Box$

\bigskip
Suppose now that the coset space $G/H$ is {\sl reductive}, i.e., $\gom$ can be
chosen such that $\Ad(h)\gom\subset \gom$ for all $h\in H$.
{}From now on assume that $\gom$ is chosen in this way.
Let
\begin{eqnarray}
\DD_H(G)&:=&\{D\in\DD(G)\mid \Ad(h)D=D\quad \mbox{for all} \ h\in H\},\\
I(\gom)&:=& \{P\in S(\gom)\mid \Ad(h)P=P\quad \mbox{for all} \ h\in H\}.
\end{eqnarray}
Then
\[\lambda(S(\gom))\cap\DD_{H,\chi,{\rm mod}}(G)=
\lambda(I(\gom))\subset\DD_H(G).\]
Hence (2.13) becomes
\begin{equation}
\DD_{H,\chi,{\rm mod}}(G)=\DD(G)\goh^\chi\oplus\lambda(I(\gom)).
\end{equation}
We obtain from Theorems 2.6 and 2.9:
\begin{thm}
Let $G/H$ be reductive.
Then:
\begin{itemize}
\item[\rm (a)] $\DD_H(G)$ is a subalgebra of $\DD(G)$.
\item[\rm (b)] $\DD(G)\goh^\chi\cap \DD_H(G)$ is a two-sided ideal in
$\DD_H(G)$.
\item[\rm (c)] There is  a direct sum decomposition
\begin{equation}
\DD_H(G)=\DD(G)\goh^\chi\cap\DD_H(G)\oplus\lambda(I(\gom)).
\end{equation}
\item[\rm (d)] Define the mappings $A, B$ and $C$ ($C$only if $\chi$ extends
to a character on $G$) by
\begin{eqnarray*}
u&\stackrel{A}{\longmapsto}& u({\rm mod}\,\DD(G)\goh^\chi\cap\DD_H(G))
\stackrel{B}{\longmapsto}
u\Big|_{C^\infty_{H,\chi}(G)}\stackrel{C}{\longmapsto} D_u\colon \\
\lambda(I(\gom))&\stackrel{A}{\longrightarrow}&
\DD_H(G)/(\DD(G)\goh^\chi\cap\DD_H(G))\stackrel{B}{\longrightarrow}
\DD_H(G)\Big|_{C^\infty_{H,\chi}(G)}\stackrel{C}{\longrightarrow}\DD(G/H).
\end{eqnarray*}
Then $A$ is a linear bijection and $B$ and $C$ are algebra isomorphisms onto.
\end{itemize}
\end{thm}

\medskip
The case $\chi\equiv 1$ of Theorem 2.10 can be found in HELGASON
[4, Cor. X.2.6 and Theor. X.2.7].
See DUFLO [1] for the general case.
\section{Application to $\DD(G/N)$ and $\DD(G/MN)$}

Let $G$ be a connected noncompact real semisimple Lie group.
We remember some of the structure theory of $G$ (cf.\ [3. Ch.VI]):

\medskip\noindent
$\gog_0$ : Lie algebra of $G$.\\
$\gog$ :  complexification of $\gog_0$.\\
$\theta$ : Cartan involution of $\gog_0$, extended to automorphism of $\gog$.\\
$\gog_0=\gok_0+\gop_0$: corresponding Cartan decomposition of $\gog_0$.\\
$\goh_{\gop_0}$: mamximal abelian subspace of $\gop_0$, $A$
the coresponding analytic subgroup.\\
$\goh_0$ : maximal abelian subalgebra of $\gog_0$ extending $\goh_{\gop_0}$.\\
$\goh_{\gok_0}:=\goh_0\cap\gok_0$,\quad $\goh_\gok$ its complexification\\
$\goh$ : complexification of $\goh_0$; this is a Cartan subalgebra of $\gog$.\\
$\Delta$ : set of roots of $\gog$ with respect to $\goh$; the roots are real on
$i\goh_{\gok_0}+\goh_{\gop_0}$.

\medskip\noindent
Introduce compatible orderings on $\goh^*_{\gop_0}$ and
$(i\goh_{\gok_0}+\goh_{\gop_0})^*$.

\medskip\noindent
$\Delta^+$ : set of positive roots.\\
$P_+$ : set of positive roots not vanishing on $\goh_{\gop_0}$.\\
$P_-$ : set of positive roots vanishing on $\goh_{\gop_0}$.\\
$\gog^{\alpha}$ : root space in $\gog$ of $\alpha\in\Delta$.\\
$\gon$ : $\sum_{\alpha\in P_+}\, \gog^\alpha$.\\
$\gon_0 :=\gon\cap \gog_0$.\\
$N$ : analytic subgroup of $G$ corresponding to $\gon_0$.\\
$M$ : centralizer of $\goh_{\gop_0}$ in $G$, $M_0$ its identity component.\\
$\gom_0$ : Lie algebra of $M$.\\
$\gom$ : complexification of $\gom_0$; then
\begin{equation}
\gom=\goh_\gok+\sum_{\alpha\in P_-} (\gog^\alpha+\gog^{-\alpha}).
\end{equation}
\begin{prop}
The coset spaces $G/MN$ and $G/N$ are not reductive.
\end{prop}
\noindent
\Proof 
Suppose that $G/MN$ is reductive.
Then there is an $\ad_\gog(\gom+\gon)$-invariant subspace
$\gor$ of $\gog$ complementary to
$\gom+\gon$.
Let $\alpha\in P_+$ and let $X$ be a nonzero element of $\gog^\alpha$.
For $H\in \goh$ write $H=W_H+Y_H+Z_H$ with $W_H\in \gor$, $Y_H\in \gom$,
$Z_H\in \gon$.
Then, for each $H\in \goh$:
\[\alpha(H)X=[W_H+Y_H+Z_H,X]\]
so
\[\alpha(H)X-[Y_H,X]-[Z_H,X]=[W_H,X]\in \gor\cap(\gom+\gon),\]
so
\[[Y_H,X]+[Z_H,X]=\alpha(H)X.\]
It follows from (3.1) that
\[[Y_H,X]+[Z_H,X]\in\sum_{\beta\in\Delta\atop \beta\neq\alpha}\, \gog^\beta.\]
Hence $\alpha(H)X=0$ for all $H\in h$, so $\alpha=0$.
This is a contradiction.

In the case $G/N$ the proof is almost the same: take $\gor$
$\ad_\gog(\gon)$-invariant and
complementary to $\gon$ and $Y_H=0$. \hfill $\Box$

\bigskip
HELGASON [5, p. 676] states without proof that $G/MN$ is not in general
reductive.

\smallskip
Let $\gol_0$ be the orthoplement of $\gom_0$ in $\gok_0$ with respect to the
Killing
form on $\gog_0$.
In order to apply Proposition 2.8 and Theorem 2.9 to $\DD(G/MN)$ and
$\DD(G/N)$ we take $\gol_0+\goh_{\gop_0}$ respectively
$\gok_0+\goh_{\gop_0}$ as
complementary subspaces of $\gom_0+\gon_0$ respectively $\gon_0$ in $\gog_0$.
Now we have
\begin{eqnarray}
I_{\rm mod}(\gol_0+\goh_{\gop_0})&=&S(\goh_{\gop_0}),\\
I_{\rm mod}(\gok_0+\goh_{\gop_0})&=&S(\gom_0+\goh_{\gop_0}).
\end{eqnarray}
(3.2) is proved in HELGASON [5. Lemma 4.2] and by only slight modifications in
this proof, (3.3) is obtained.
It follows from Lemma 2.5 that
\[\lambda(S(\goh_{\gop_0}))\subset\DD_{MN,1,{\rm mod}}(G)\]
and
\[\lambda(S(\gom_0+\goh_{\gop_0}))\subset\DD_{N,1,{\rm mod}}(G),\]
since $M$  centralizes $\goh_{\gop_0}$ and $\gom_0+\goh_{\gop_0}$
normalizes $\gon_0$.
Consider $\DD(A)$ and $\DD(M_0A)$ as subalgebras of $\DD(G)$.
Then $\DD(A)\subset\DD_{MN,1,{\rm mod}}(G)$ and
$\DD(M_0A)\subset\DD_{N,1,{\rm mod}}(G)$.
It follows by application of Proposition 2.8 and Theorem 2.9 that:
\begin{thm}
The mapping $u\mapsto D_u$ (cf. (2.17)) is an algebra isomorphism from $\DD(A)$
onto $\DD(G/MN)$ and from $\DD(M_0A)$ onto $\DD(G/N)$.
In particular, $\DD(G/MN)$ is a commutative algebra.
\end{thm}

\medskip
The statements about $\DD(G/MN)$ are in HELGASON [5, Theorem 4.1].
FARAUT [2, p.~393] observes that Helgason's result can be extended to the
context
of pseudo-Riemannian symmetric spaces.

A special case of Theorem 6.2 can be formulated in the situation that $G$ is a
connected complex semisimple Lie group.
Let $\gog$ be its (complex) Lie algebra and put:

\medskip\noindent
$\gou$ : compact real form of $\gog$.\\
$\goa$ : maximal abelian subalgebra of $\gou$.\\
$\goh:= \goa+i\goa$; this is  Cartan subalgebra of $\gog$.\\
$\Delta$ : set of roots of $\gog$ with respect to $\goh$.\\
$\Delta^+$ : set of positive roots with respect to some ordering.\\
$\gog^\alpha$ : root space of $\alpha\in\Delta$.\\
$\gon:=\sum_{\alpha\in\Delta^+}\gog^\alpha$,\quad $N$ the corresponding
analytic subgroup.\\
$\gog^\RR:=\gog$ considered as real Lie algebra.\\
$\goh^\RR:=\goh$ considered as real subalgebra.

\medskip\noindent
Then $\gog^\RR=\gou+i\goa+\gon$ is an Iwasawa decomposition for
$\gog^\RR$ (cf.\ [4, Theorem
VI.6.3]) and $\goa$ is
the centralizer of $i\goa$ in $\gou$.
Hence we obtain from Theorem~3.2:
\begin{thm}\quad
The mapping $P\mapsto D_{\lambda(P)}$ is an algebra isomorphism from
$S(\goh^\RR)$
onto $\DD(G/N)$.
In particular, $\DD(G/N)$ is commutative.
\end{thm}

\medskip
This theorem was proved by HELGASON [6, Lemma 3.3] without use of Theorem 3.2.

\section{Application to isotropic spaces}

We preserve the notation and conventions of Section 2.
First we prove an extension of [4, Cor. X.2.8] to the case that $G/H$ is
not necessarily
reductive.
In the following, $A$ and $B$ are as in Theorem 2.6(d).
\begin{lem}
If the algebra $I_{\rm mod}(\gom)$ is generated by $P_1,\dots,P_l$ and if
there are
$Q_1,\dots,Q_l\in S_m$ such that ${\rm degree}(P_i-Q_i)<{\rm degree}\,
P_i$ and
$\lambda(Q_i)\in\DD_{H,\chi,{\rm mod}}(G)$ then the algebra
\[\DD_{H,\chi,{\rm mod}}\Big|_{C^\infty_{H,\chi}(G)}\]
is generated by $BA\lambda(Q_1),\dots,BA\lambda(Q_l)$.
\end{lem}

\medskip\noindent
\Proof 
We prove by complete induction with respect to $m$ that, for each $P\in
S_m(\gom)$ with
$\lambda(P)\in \DD_{H,\chi,{\rm mod}}(G), \ BA\lambda(P)$ depends
polynomially on
$BA\lambda(Q_1),\dots, BA\lambda(Q_l)$.
In view of Theorem 2.6 this will prove the lemma.
Suppose the above property holds up to $m-1$.
Let $P\in S_m(\gom)$ such that $\lambda(P)\in\DD_{H,\chi,{\rm mod}}(G)$.
By using Lemma 2.7 we find that $P=\Pi(P_1,\dots, P_l) \pmod
{S_{m-1}(\gom)}$ for
some polynomial $\Pi$ in $l$ indeterminates.
Hence, $P=\Pi(Q_1,\dots, Q_l)  \pmod {S_{m-1}(\gom)}$,
\begin{eqnarray*}
\lambda(P)&=&\lambda(\Pi(Q_1,\dots,Q_l)) \pmod{\lambda(S_{m-1}(\gom))}\\
&=&\Pi(\lambda(Q_1),\dots,\lambda(Q_l)) \pmod{\lambda(S_{m-1}(\gog))},\\
\lambda(P)&-&\Pi(\lambda(Q_1),\dots,\lambda(Q_l))\in\lambda(S_{m-1}(\gog))
\cap\DD_{H,\chi,{\rm mod}}(G).
\end{eqnarray*}
By Theorem 2.6 and formula (2.10) we have
\[BA\lambda(P)-\Pi(BA\lambda(Q_1),\dots, BA\lambda(Q_l))=BA\lambda(P')\]
for some $P'\in S_{m-1}(\gom)$ such that $\lambda(P')\in\DD_{H,\chi,{\rm
mod}}(G)$.
Now apply the induction hypothesis.\\ \mbox{}\hfill $\Box$

\medskip
Let $\tau$ denote the action of $G$ on $G/H$.
Its differential $d\tau$ yields an action of $H$ on the tangent space
$(G/H)_0$ to $G/H$ at $0$.
\begin{thm}
Suppose there is a nondegenerate $d\tau(H)$-invariant bilinear form
$\langle\cdot,\cdot\rangle$ on $(G/H)_0$ of signature $(r_1,r_2) \
(r_1+r_2=r, \ r_1\ge r_2)$
such that, for each $\lambda>0, \ d\tau(H)$ acts transitively on $\{X\in
(G/H)_0\mid\langle
X,\, X\rangle=\lambda\}$ (or on the connected components of these
hyperbolas if $r_1=r_2=1$).
Let $\Delta$ be the Laplace-Beltrami operator on $G/H$ coresponding to the
$G$-invariant
pseudo-Riemannian structure on $G/H$ associated with
$\langle\cdot,\cdot\rangle$.
Then the algebra $\DD(G/H)$ is generated by $\Delta$, and hence commutative.
\end{thm}

\medskip\noindent
\Proof 
Choose a complementary subspace $\gom$ to $\goh$ in $\gog$.
The mapping $d\pi$ identifies the $H$-spaces $\gom$ (under $\sigma\circ
\Ad_G(H)$) and $(G/H)_0$
(under  $d\tau(H)$) with each other.
Transplant the form $\langle\cdot,\cdot\rangle$ to $\gom$ and choose an
orthonormal basis
$X_1,\dots,X_r$ of $\gom$: $\langle X_i,\, X_j\rangle=
\varepsilon_i\delta_{ij}$, \
$\varepsilon_i=1$ or $-1$ for $i\le r_1$ or $>r_1$, respectively.
Then the algebra $I_{\rm mod}(\gom)$ is generated by
$\sum^r_{u=1}\varepsilon_iX_i^2$.
It follows from the proof of Theorem 2.9 that $\Delta=D_{\lambda(P)}$ with
$P\in S(\gom)$ of
degree 2 such that $\lambda(P)\in\DD_{H,1,{\rm mod}}(G)$.
Thus, by Lemma 2.7, we get
\[P=c\sum_{i=1}^r\, \varepsilon_iX_i^2\pmod{S_1(\gom)}\]
with $c\neq 0$.
Now apply Lemma 4.1 and Theorem 2.9.\hfill $\Box$

\bigskip
Theorem 4.2 extends [4, Prop.\ X.2.10], where the case is considered that
$G/H$ is a Riemannian
symmetric space of rank $1$.
A pseudo-Riemannian manifold  $M$ is called {\sl isotropic} if for each
$x\in M$ and for
tangent vectors $X,Y\neq 0$ at $x$ with $\langle X,\, X\rangle=\langle Y,\,
Y\rangle$ there
is an isometry of $M$ fixing $x$ which sends $X$ to $Y$.
Connected isotropic spaces can be written as homogeneous spaces $G/H$
satisfying the
conditions of Theorem 4.2 with $G$ being the full isometry group (cf. WOLF
[8, Lemma 11.6.6]).
It follows from Wolf's classification [8, Theorem 12.4.5] that such spaces
are symmetric and
reductive.
However, our proof of Theorem 4.2 does not use this fact.

\vskip 0.5truecm
\noindent
{\obeylines\parindent 9.7truecm
\addressfont
{\sladdressfont present address:}
Korteweg-de Vries Institute for Mathematics
Universiteit van Amsterdam
Plantage Muidergracht 24
1018 TV Amsterdam, The Netherlands
email: {\ttaddressfont thk@science.uva.nl}}

\end{document}